\documentclass[MSNbibl,number,citesort,dvips]{arxbj}
\usepackage{upgreek}

\aid{0}
\volume{19}
\issue{2}
\pubyear{2013}
\firstpage{387}
\lastpage{408}
\doi{10.3150/11-BEJ405}

\makeatletter
\newcommand\R{\mathbb{R}}

\newcommand\N{\mathbb{N}}
\newcommand\W{\mathbb{N}_0}
\newcommand\HE{\operatorname{HE}}
\newtheorem{theorem}{Theorem}

\newcommand\corr{\operatorname{corr}}

\newcommand\bx{\mathbf{x}}
\newcommand\by{\mathbf{y}}

\newcommand\bz{\mathbf{z}}

\newcommand\ze{\mathbf{0}}

\newcommand\bnu{\bolds{\nu}}

\newcommand\bom{\bolds{\omega}}
\newcommand\bdom{\,\mathrm{d}\bolds{\omega}}
\newcommand\btau{\bolds{\tau}}
\newcommand\bdtau{\,\mathrm{d}\bolds{\tau}}
\renewcommand\l{\lambda}

\newcommand\Var{\operatorname{Var}}
\makeatother

\begin{document}
\begin{frontmatter}

\title{On a class of space--time intrinsic random functions}
\runtitle{Space--time intrinsic random functions}

\begin{aug}
\author{\fnms{Michael L.} \snm{Stein}\corref{}\ead[label=e1]{stein@galton.uchicago.edu}}
\runauthor{M.L. Stein} 
\address{Department of Statistics, University of Chicago,
Chicago, IL 60637, USA.\\
\printead{e1}}
\end{aug}

\received{\smonth{3} \syear{2011}}
\revised{\smonth{8} \syear{2011}}

%
\begin{abstract}
Power law generalized covariance functions provide a simple model
for describing the local behavior of an isotropic random field.
This work seeks to extend this class of covariance functions to
spatial-temporal processes for which the degree of smoothness in
space and in time may differ while maintaining other desirable
properties for the covariance functions, including the availability
of explicit convergent and asymptotic series expansions.
\end{abstract}

%

\end{frontmatter}

\section{Introduction}\label{sec1}
Intrinsic random functions~\cite{chiles,matheron} provide a popular
class of
models for spatial processes.
These non-stationary random processes are specified by their generalized
covariance functions, which determine the variances of certain linear
combinations of the process (see Section~\ref{sec2.1} for details).
A particularly simple and, therefore, useful class of generalized covariance
functions is the power law class, for which the ``covariance'' between two
observations is proportional to the Euclidean
distance between the points raised to some power (unless the power is
an even
integer).
More specifically, for $x\ge0$ indicating interpoint distance,
$\zeta>0$ and $\N$ the set of positive integers,
the function $\gamma_\zeta$
%
\begin{equation}\label{gamma.def}
\gamma_\zeta(x) =
\cases{ \Gamma(-\zeta) x^{2\zeta},
& \quad $\zeta\notin\N$,\vspace*{2pt}\cr
\displaystyle\frac{2(-1)^{\zeta+1}}{\zeta!}x^{2\zeta}\log x, &\quad$\zeta\in\N$}
\end{equation}
gives a valid generalized covariance function in any number of dimensions.
Despite its simplicity, this class of models has the important virtue
of admitting a broad range of local behaviors for the process,
which is critical when, for example, considering properties of spatial
interpolants~\cite{stein-book}.
Specifically,
the larger the value of $\zeta$, the smoother the process, so that,
for example, the process is $m$
times mean square differentiable in any direction if and only if $\zeta>m$.
A generalized covariance function can be written
as the Fourier transform of a positive measure; although, in contrast
to the stationary setting, the measure might not have finite total mass.
In particular, in~$d$ dimensions,
the measure corresponding to $\gamma_\zeta$ has density
with respect to Lebesgue measure
proportional to $|\bom|^{-2\zeta-d}$, where $\bom\in\R^d$ is the spatial
frequency.

The goal of this paper is to find a good extension of $\gamma_\zeta$
to the space--time setting.
The problem is made difficult by what I will mean by ``good.''
The first requirement is that the class of models includes members
allowing any degree of smoothness in space and any (possibly different) degree
of smoothness in time.
Specifically, for
$G(\bx,t)$ a generalized covariance function of spatial lag $\bx$ and
temporal lag $t$, and any positive $\zeta_1,\zeta_2,C_1$ and $C_2$,
the class of generalized covariance functions should
include a member that satisfies
%
\begin{equation}
G(\bx,0) = C_1\gamma_{\zeta_1}(|\bx|)
\label{coord1}
\end{equation}
and
%
\begin{equation}
G(0,t) = C_2\gamma_{\zeta_2}(|t|)
\label{coord2}
\end{equation}
or, failing that, that (\ref{coord1}) holds asymptotically as $\bx\to\ze$,
and (\ref{coord2}) holds asymptotically as $t\to0$.

The second requirement is that $G$ be smoother away from the origin than
it is at the origin.
This requirement is needed to avoid the kinds of anomalies described in
\cite{stein2005}
for covariance functions that are not smoother
away from the origin.
As a simple example of this kind of anomaly, consider the covariance
function on $\R\times\R$ given by $K(x,t) = \exp(-|x|-|t|)$.
Write corr for correlation, and
define $\rho(x,t) = \lim_{\varepsilon\downarrow0}\operatorname{corr}\{Z(0,\varepsilon)-
Z(0,0), Z(x,t+\varepsilon)-Z(x,t)\}$.
Then for $x\ne0$, straightforward calculations show $\rho(x,t)=0$ for
$t\ne
0$ and $\rho(x,0)=\mathrm{e}^{-|x|}$.
The discontinuity in this limiting correlation is due to the fact that
$K(x,t)$ has a similar discontinuity in its first derivative in the $t$
direction everywhere along the $x$ axis as it does at the origin.
I consider such a discontinuity in $\rho$
as unrealistic for most natural space--time processes.
In particular, Stein and Handcock
\cite{stein-handcock} and Stein~\cite{stein2011} give examples showing
how this
lack of continuity away from the origin can lead to optimal (kriging)
predictors with undesirable properties.
All separable space--time covariance functions, that is, those that
factor into a function of space and a function of time such as $\mathrm{e}^{-|x|-|t|}$,
have a similar
problem unless the process is infinitely differentiable~\cite{stein2005}, page~311.
Furthermore, many non-separable space--time covariance functions proposed
in the literature share this problem~\cite{stein2005}, pages~311--312.
For a space--time covariance function with different degrees of
smoothness in
space and time, that is, satisfying (\ref{coord1}) and (\ref{coord2}) with
$\zeta_1\ne\zeta_2$, it is not so clear what one should mean by the function
being smoother away from the origin than at the origin.
This issue is addressed in Section~\ref{sec2.4}.

The smoothness of a covariance function away from the origin is
closely related to regularity properties of the corresponding spectral
density at high frequencies~\cite{stein2005}.
Indeed, Stein~\cite{stein2011} argues that the spectral domain
provides a more
natural approach for considering the appropriateness of various models
for space--time covariance functions.
Specifically,~\cite{stein2011} gives the following condition as
a plausible requirement for the spectral density $f(\bom)$
of a natural process in space or space--time:
for every $R<\infty$,
%
\begin{equation}
\lim_{\bom\to\infty} \sup_{|\bnu|<R} \biggl| \frac{f(\bom+\bnu)}{f(\bom)}
-1\biggr| = 0.
\label{f-cond}
\end{equation}
That is, $f$ changes slowly (on a relative scale) at high
frequencies.
This condition excludes, for example, separable space--time models.

The third and final requirement for the generalized covariance
functions is that they can be computed accurately and efficiently via,
for example, series expansions,
to allow them to be applied routinely to large space--time
datasets.
In particular, representations of the function as an integral
will not be considered an adequate solution to the problem.
As best as I am aware, no existing class of generalized covariance functions
satisfies all three of these requirements.

It will be convenient to avoid explicitly distinguishing between
space and time and consider processes $Z(\bx,\by)$ on $\R^{d_1}\times
\R^{d_2}$ for positive integers $d_1$ and $d_2$.
Stein~\cite{stein2005}
proposed the following class of spectral densities (with a different
parameterization) as a flexible
parametric model for stationary space--time processes:
%
\begin{equation}
q(\btau,\bom)=\biggl\{ \biggl(\frac{\beta_1^2+|\btau|^2}{\sigma_1^2}
\biggr)^{\alpha_1}
+ \biggl(\frac{\beta_2^2+|\bom|^2}{\sigma_2^2}
\biggr)^{\alpha_2}\biggr\} ^{-\nu}
\label{model1}
\end{equation}
for $\sigma_1,\sigma_2$ positive, $\beta_1^2+ \beta_2^2>0$
and $2\nu>d_1/\alpha_1+d_2/\alpha_2$; this last condition
being necessary and sufficient (given the positivity constraints
on the other parameters) for $q$ to be integrable.
The parameters $\beta_1$ and $\beta_2$ are inverse
range parameters, $\sigma_1$ and $\sigma_2$
are scale parameters and $\alpha_1,\alpha_2$ and $\nu$ together
describe the smoothness of the process in $\bx$ and $\by$.
More specifically, the process is $p$ times mean square differentiable in
each component of $\bx$ if and only if $2\nu>(d_1+2p)/\alpha
_1+d_2/\alpha_2$
and, similarly,
is $p$ times mean square differentiable in components
of $\by$ if and only if $2\nu>d_1/\alpha_1+(d_2+2p)/\alpha_2$
\cite{stein2005}.
Furthermore, when $\alpha_1$ and $\alpha_2$ are integers, the resulting
covariance function is infinitely differentiable away from the origin
(\cite{stein2005}, Proposition~4, although this result also follows from
\cite{shubin}, Theorem 1.1).
All models in the class (\ref{model1}) satisfy the spectral condition
(\ref{f-cond}).
Porcu~\cite{porcu11} describes more general approaches to obtaining
valid spectral
densities that could be useful in the space--time context.

An obstacle to using the class of models (\ref{model1})
is the lack of an explicit expression for the corresponding covariance
functions.
Except when $\alpha_1=\alpha_2=1$ and some special cases with
$\nu$ an integer, $\alpha_1=1$ and $\alpha_2=2$, I am unaware of any
cases for which an explicit expression has been written down~\cite{stein2005}.
For rational spectral densities, which includes the model
(\ref{model1}) when $\alpha_1,\alpha_2$ and~$\nu$ are all integers as
a special case, the covariance function can be expressed as the solution
of a certain set of equations that reduce to
a partial differential equation when the rational function
is just a reciprocal of a polynomial~\cite{ramm03,ramm05}.
However, writing down an explicit solution
for all integers $\alpha_1,\alpha_2$ and $\nu$ is not a simple task.
Ma~\cite{ma2007} gives explicit expressions for the covariance functions
of space--time processes with rational spectral densities in some
limited special cases that do not include any instances of (\ref{model1}).
For all of the cases treated
in~\cite{ma2007}, the covariance functions are not smoother away from
the origin than at the origin in the sense defined in Section~\ref{sec2.3},
and the spectral densities do not satisfy (\ref{f-cond}).
Even if one had explicit expressions for
the covariance functions of all rational spectral
densities, these covariance functions
(asymptotically in a neighborhood of the origin)
satisfy (\ref{coord1}) and (\ref{coord2}) only for a countable
nowhere dense set of $(\zeta_1,\zeta_2)$ values.
Kelbert, Leonenko and Ruiz-Medina~\cite{kelbert} consider the case $d_1=d_2=1,
\alpha_1=1,\alpha_2=2$ and $\beta_1=0$ as a stochastic fractional heat equation
in some detail, but they only give integral representations for the covariance
functions.
They consider additionally setting $\beta_2=0$
and note that the resulting random
field has a self-similarity property when $\nu\in(\frac{3}{4},\frac{7}{4}
)$.
Christakos~\cite{christakos00}, page~225,
mentions the case with $d_1=1,\alpha_1=\nu=1$
and $\beta_1=0$ as a possible model for spatial-temporal processes, but derives
no results for these models.

The parameters $\beta_1$ and $\beta_2$
are range parameters that do not affect the
local behavior of the process, so consider setting
$\beta_1=\beta_2=0$ in (\ref{model1}), yielding
%
\begin{equation}
f(\btau,\bom) = \biggl\{\biggl(\frac{|\btau|}{\sigma_1}\biggr)^{2\alpha
_1} +
\biggl(\frac{|\bom|}{\sigma_2}\biggr)^{2\alpha_2}\biggr\}^{-\nu}.
\label{model2}
\end{equation}
These spectral densities are not integrable in a neighborhood
of the origin, so they correspond to generalized covariance
functions.
Section~\ref{sec2.1} gives some background on generalized covariance functions.
Theorem~\ref{Gtheorem} in
Section~\ref{sec2.2} gives a convergent power series for the generalized covariance
function corresponding to
(\ref{model2}) when $\alpha_2=1$, $\alpha_1 > 1$ and $\by\ne\ze$ and
separate explicit formulae to cover the case $\by= \ze$.
For completeness, the known result~\cite{gelfand}
for the generalized covariance
function when $\alpha_1=1$ is also given in Theorem~\ref{Gtheorem}.

Section~\ref{sec2.3} shows how the generalized covariance functions in Theorem
\ref{Gtheorem}
can, in most cases,
be written in terms of the $H$-function~\cite{kilbas},
a generalization of the generalized hypergeometric function
that is sometimes called Fox's $H$-function~\cite{fox}.
This result is used to obtain asymptotic series for these generalized
covariance functions and to motivate a
conjecture extending Theorem~\ref{Gtheorem} to the case $\alpha_1<1$.
Section~\ref{sec2.4} shows that for any positive $C_1,C_2,\zeta_1$
and $\zeta_2$, one can find a spectral density of this form
satisfying (\ref{coord1}) and (\ref{coord2}).
Furthermore, in a sense made precise in Section~\ref{sec2.4}, the resulting
generalized covariance function is shown to be
smoother away from the origin than at the origin.
Finally, these results are used to show that
at least some of these covariance
functions avoid what~\cite{kent} calls the ``dimple'' that occurs
in some proposed space--time covariance functions, which is a lack of
monotonicity in the correlation structure that might often be viewed
as unnatural.

Section~\ref{sec3} discusses
some limitations and possible extensions of the generalized covariance
functions considered in Section~\ref{sec2}.
This section also touches on some of the
difficulties in using the series expansions to compute these functions
quickly and accurately.
Section~\ref{appa1} gives a proof of Theorem~\ref{Gtheorem}, and Section~\ref{appa2} collects some
needed material on $H$-functions.

There is a substantial recent literature on the development of space--time
covariance functions with explicit representations in terms of well-known
special functions.
Some references include
\cite{christakos00,cressie,gneiting,gregori,jones,kolovos,ma2003,ma2005,ma2007,ma2008,porcu,porcu08,schlather,stein2005,zastavnyi}.
Although these works consider a broad range of models for space--time
covariance functions, none of them give a class of covariance functions
meeting the criteria set forth in this section of
(locally) satisfying
(\ref{coord1}) and (\ref{coord2}) for all positive $\zeta_1$ and $\zeta_2$
as well as being smoother away from the origin than at the origin.
As noted in~\cite{stein2005}, perhaps~\cite{ma2003}
comes closest to this goal,
in that this paper gives a class of models with $d_2=1$
including elements satisfying,
asymptotically in a neighborhood of the origin,
(\ref{coord1}) with $\zeta_1=\frac{1}{2}$ and (\ref{coord2})
with $0<\zeta_2<\frac{1}{4}$, and
the covariance functions are infinitely differentiable away from the origin.

\section{Theoretical results}\label{sec2}

\subsection{Generalized space--time covariance functions}\label{sec2.1}
Intrinsic random functions
and generalized covariance functions have been a standard tool in
geostatistics since Matheron's pioneering work~\cite{matheron}.
These processes are nearly stationary in the sense that variances of
some class of linear combinations of the process are translationally
invariant.
Specifically,
for a random field $Z$ on $\R^d$ (so that $d=d_1+d_2$ is the total number
of dimensions for space--time processes),
call $\sum_{j=1}^n \l_j Z(\bz_j)$
an authorized linear combination of order $k$, or ALC-$k$, if
$\sum_{j=1}^n \l_jP(\bz_j) = 0$ for every polynomial $P$ of
order at most $k$.
A function $G$ on $\R^d$ is called a generalized covariance function of order
$k$, or GC-$k$, if, for every ALC-$k$,
$\sum_{\ell,j=1}^n \l_\ell\l_jG(\bz_\ell-\bz_j)\ge0$.
A process $Z$ for which $\Var\{\sum_{j=1}^n \l_j Z(\bz_j)\}
= \sum_{\ell,j=1}^n
\l_\ell\l_jG(\bz_\ell-\bz_j)\ge0$
and $E\{\sum_{j=1}^n \l_j Z(\bz_j)\} = 0$
for every ALC-$k$ is said to be an
intrinsic random function of order $k$, or IRF-$k$,
with $G$ as its GC-$k$.
GC-$k$s are not unique; if $G$ is a GC-$k$ for $Z$, then so is $G$ plus
any even polynomial of degree $2k$ in $\bx$.
A function $f$ is the spectral density of a GC-$k$ $G$ if
\[
\sum_{\ell,j=1}^n \l_\ell\l_jG(\bz_\ell-\bz_j)
= \int_{\R^d} \Biggl|\sum_{j=1}^n \l_j \mathrm{e}^{\mathrm{i}\bom'\bz_j}\Biggr|^2 f(\bom)
 \bdom
\]
for every ALC-$k$ $\sum_{j=1}^n \l_j Z(\bz_j)$.
A nonnegative even function $f$ is the spectral density for a real-valued
GC-$k$ if and only if $f(\bom)|\bom|^{2k+2}/
(1+|\bom|^{2k+2})$ is integrable~\cite{matheron}.
Of course, if $Z$ is an IRF-$k$, it is also an IRF-$k'$ for all integers
$k'\ge k$.
Let $\lfloor x\rfloor$ indicate
the greatest integer less than or equal to $x$.
For the spectral density (\ref{model2}) and $i=1,2$, define
$k_i = \lfloor\alpha_i\{\nu-d_1/(2\alpha_1) -d_2/(2\alpha_2)\}\rfloor$
and $k_0 = \max(k_1,k_2)$.
Straightforward calculus shows that $f$ in (\ref{model2}) satisfies
the integrability condition for an IRF-$k$ if and only if $k \ge k_0$.
Christakos~\cite{christakos}
considers an extension of the notion of IRFs to the space--time
setting in which one essentially allows a different degree of
differencing in
space and in time, but this concept is not used here.

\subsection{Main theorem}\label{sec2.2}

This section gives a series expansion
for the generalized covariance function corresponding to~(\ref{model2}) when $\alpha_2=1$, $\alpha_1>1$ and $\sigma_1=\sigma_2=1$.
Extending the result to other positive values of $\sigma_1$ and $\sigma
_2$ is
trivial.\vadjust{\goodbreak}
Define the function $\mathcal{M}_\nu(y) = y^\nu\mathcal{K}_\nu(y)$, where $\mathcal{K}_\nu$ is a modified Bessel function.
When $\nu>0$, $\mathcal{M}_\nu$ is often called
the Mat\'ern covariance function with smoothness parameter $\nu$.
Set $\theta= \nu-\frac{1}{2}d_2$ and
$\theta'=\theta-d_1/(2\alpha_1)$
so that $k_0 = \lfloor\alpha_1\theta' \rfloor$ is the (minimal)
order of the IRF corresponding to this spectral density.
As noted in~\cite{stein2005} for the more general model (\ref{model1}),
when $\alpha_2=1$, the Fourier
transform with respect to $\bom$ can be carried out explicitly.
Specifically, for $\btau\ne\ze$,
%
\begin{equation}
\int_{\R^{d_2}} ( |\btau|^{2\alpha_1}+
|\bom|^2)^{-\nu} \mathrm{e}^{\mathrm{i}\bom' \by} \bdom
= \frac{\uppi^{d_2/2} \mathcal{M}_{\theta}
(|\btau|^{\alpha_1}|\by|)}
{2^{\theta-1}\Gamma(\nu) |\btau|^{2\alpha_1\theta}}.
\label{partint}
\end{equation}
Define $r_{\ell j} = |\bx_\ell-\bx_j|$, $s_{\ell j} =
|\by_\ell-\by_j|$ and, for $t>0$,
%
\begin{equation}
\Lambda_d(t) = 2^{(d-2)/2}\Gamma\biggl(\frac{1}{2}d\biggr)
t^{-(d-2)/2} J_{(d-2)/2}(t) = \sum_{m=0}^\infty
\frac{(-({1}/{4})t^2)^m}{m!({d}/{2})_m},
\label{hankel}
\end{equation}
where, for any real $a$ and
positive integer $j$, $(a)_j = a(a+1)\cdots(a+j-1)$ and
$(a)_0=1$.
For every ALC-$k_0$ $\sum_{j=1}^n \l_j Z(\bx_j,\by_j)$,
the GC-$k_0$ $G$ corresponding to the spectral density
$ ( |\btau|^{2\alpha_1} + |\bom|^2)^{-\nu}$ satisfies
%
\begin{eqnarray}\label{step0}
& & \sum_{\ell,j=1}^n \l_\ell\l_j G(\bx_\ell-\bx_j,\by_\ell-\by_j)
\nonumber\\
& & \quad  =
\int_{\R^{d_1}}\int_{\R^{d_2}} \sum_{\ell,j=1}^n \l_\ell\l_j
\mathrm{e}^{\mathrm{i}\btau'(\bx_\ell-\bx_j) + \mathrm{i}\bom'(\by_\ell-\by_j)}
(|\btau|^{2\alpha_1}+|\bom|^2)^{-\nu}\bdom \bdtau
\nonumber\\
& & \quad  = \frac{\uppi^{d_2/2}}{2^{\theta-1}\Gamma(\nu)}
\int_{\R^{d_1}}\sum_{\ell,j=1}^n \l_\ell\l_j
\frac{\mathcal{M}_{\theta} (|\btau|^{\alpha_1}s_{\ell j})}
{|\btau|^{2\alpha_1\theta}}
\mathrm{e}^{\mathrm{i}\btau'(\bx_\ell-\bx_j)}\bdtau
\nonumber
\\[-8pt]
\\[-8pt]
\nonumber
& & \quad  = \frac{4\uppi^{(d_1+d_2)/2}}{2^\theta\Gamma(\nu)\Gamma(
({1}/{2})d_1)}
\int_0^\infty\sum_{\ell,j=1}^n \l_\ell\l_j u^{d_1-1-2\alpha_1\theta}
\Lambda_{d_1}(r_{\ell j}u)\mathcal{M}_\theta(s_{\ell j}u^{\alpha_1})
  \,\mathrm{d}u \\
& & \quad  = \frac{4\uppi^{(d_1+d_2)/2}}{2^\theta\Gamma(\nu)\Gamma(
({1}/{2})d_1)\alpha_1} \nonumber\\
& & \qquad{} \times
\int_0^\infty\sum_{\ell,j=1}^n \l_\ell\l_j t^{d_1/\alpha_1-2\theta-1}
\Lambda_{d_1}(r_{\ell j}t^{1/\alpha_1})\mathcal{M}_\theta(s_{\ell j}t)
  \,\mathrm{d}t,\nonumber
\end{eqnarray}
where the second step uses (\ref{partint}), the third
basic results on Fourier transforms of isotropic functions~\cite{stein-book}, Section 2.10, and the last step
the change of variables $t=u^{\alpha_1}$.
What one would like to do is substitute (\ref{hankel}) into (\ref{step0})
and integrate termwise, but justifying this interchange requires
considerable care.

Define
\[
c_m(\alpha_1) = \frac{\uppi^{(d_1+d_2)/2}
\Gamma({(d_1+2m)}/{(2\alpha_1)})}
{m!\Gamma(m+({1}/{2})d_1)}
\]
and the digamma function $\psi$ by $\psi(z) = \frac{\mathrm{d}}{\mathrm{d}z}\log\Gamma(z)$.
Equations (\ref{Gres})--(\ref{G1k0.1}) are proven in Section~\ref{appa1}.
Equations (\ref{isotropic}) and (\ref{isotropic2}) are taken from
\cite{gelfand}, Chapter II, Section 3.3, equations (2) and~(11).

\begin{theorem}\label{Gtheorem}
For the spectral density $f(\btau,\bom) = (|\btau|^{2\alpha_1} + |\bom
|^2)^{-\nu}$
on $\R^{d_1}\times\R^{d_2}$ with $\alpha_1\ge1$,
$\theta'=\nu-d_1/(2\alpha_1)-d_2/2>0$ and $k_0=
\lfloor\alpha_1\theta'\rfloor$, a corresponding GC-$k_0$ is given by
$\tilde G(
|\bx|,|\by|)$ for the function $\tilde G$, defined by the following equations.
First consider $\alpha_1>1$.
For $s>0$,
%
\begin{equation}
\tilde G(r,s) =
\sum_{m=0}^\infty\frac{c_m(\alpha_1)}{\alpha_1\Gamma(\nu)}
\biggl\{-\biggl(\frac{1}{2}r\biggr)^2\biggr\}^m \gamma_{\theta'-m/\alpha_1}
\biggl(\frac{1}{2}s\biggr)
\label{Gres}
\end{equation}
with $\gamma(\cdot)$ defined by (\ref{gamma.def}).
When $\alpha_1\theta' > k_0$,
%
\begin{equation}
\tilde G(r,0) = \frac{\uppi^{(d_1+d_2)/2}\Gamma(\theta)}
{\Gamma(\nu)\Gamma(\alpha_1\theta)}
\gamma_{\alpha_1\theta'}\biggl(\frac{1}{2}r\biggr),
\label{G1k0.0}
\end{equation}
and when $\alpha_1\theta' = k_0$,
%
\begin{eqnarray}\label{G1k0.1}
\tilde G(r,0) & = & - \frac{c_{k_0}(\alpha_1)}{\Gamma(\nu_0)}
\biggl\{-\biggl(\frac{1}{2}r\biggr)^2\biggr\}^{k_0}
\biggl\{2\log\biggl(\frac{1}{2}r\biggr)+
\frac{1}{\alpha_1}\psi\biggl(\frac{2k_0+d_1}{2\alpha_1}\biggr)
\nonumber
\\[-8pt]
\\[-8pt]
\nonumber
& &\phantom{ - \frac{c_{k_0}(\alpha_1)}{\Gamma(\nu_0)}
\biggl\{-\biggl(\frac{1}{2}r\biggr)^2\biggr\}^{k_0}
\biggl\{}{}  +\frac{1}{\alpha_1}\psi(1) - \psi\biggl(k_0+\frac
{1}{2}d_1\biggr)
- \psi(k_0+1) \biggr\}.
\end{eqnarray}
Finally,
when $\alpha_1=1$ and $\theta'$ is not an integer,
%
\begin{equation}
\tilde G(r,s) =
-\frac{\uppi^{(d_1+d_2+2)/2}}{\sin(\uppi\theta')\Gamma(\theta'+1)
\Gamma(\nu)2^{2\theta'}} (r^2+s^2)^{\theta'}
\label{isotropic},
\end{equation}
and when $\alpha_1=1$ and $\theta'$ is an integer,
%
\begin{equation}
\tilde G(r,s) =
\frac{(-1)^{k+1}\uppi^{(d_1+d_2)/2}}{\theta'!
\Gamma(\nu)2^{2\theta'}} (r^2+s^2)^{\theta'}\log(r^2+s^2).
\label{isotropic2}
\end{equation}
\end{theorem}

\subsection{$H$-functions}\label{sec2.3}

The function $\tilde G$ can, in most cases, be written in terms of
$H$-functions~\cite{kilbas,mathai}.
Section~\ref{appa2} gives the definition of $H$-functions as a contour integral
and some other needed information about the functions.
If
%
\begin{equation}
\ell+\frac{m}{\alpha_1}\neq\theta'
\qquad\mbox{for all whole numbers } \ell\mbox{ and } m,
\label{Hcondition}
\end{equation}
then (\ref{poles}) in Section~\ref{appa2} is satisfied, and
\[
H^{2,1}_{2,2} \left( z \left|
\matrix{(1,1), \biggl(\displaystyle\frac{d_1}{2},1\biggr) \vspace*{2pt}\cr
\biggl(\displaystyle\frac{d_1}{2\alpha_1},\displaystyle\frac
{1}{\alpha_1}\biggr),\biggl(-\theta',\frac{1}{\alpha_1}\biggr)}\right.\right)
\]
is well defined.
Since, by (\ref{Delta}), $\Delta=\frac{2}{\alpha_1}-2$,
which is negative for $\alpha_1>1$,
\cite{kilbas}, Theorem~1.4, applies, yielding, for $z\ne0$,
%
\begin{equation}
H_{2,2}^{2,1}\left(z\left|
\matrix{(1,1), \biggl(\displaystyle\frac{d_1}{2},1\biggr)\vspace*{2pt}\cr
\biggl(\displaystyle\frac{d_1}{2\alpha_1},\displaystyle\frac
{1}{\alpha_1}\biggr),\biggl(-\theta',\frac{1}{\alpha_1}\biggr)}\right.\right)
= \sum_{k=0}^\infty h_{1k}z^{-\alpha_1 k},
\label{Hseries}
\end{equation}
where
\[
h_{1k} = {(-1)^k\alpha_1\Gamma\biggl(\frac{d_1+2k}{2\alpha_1}\biggr)
\Gamma\biggl(-\theta'+\frac{k}{\alpha_1}\biggr)}\Big/{\biggl(k!\Gamma\biggl(\frac{d_1}{2}+k
\biggr)\biggr)}.
\]
Comparing this result to (\ref{Gres}) yields (assuming (\ref{Hcondition}))
%
\begin{equation}
\tilde G(r,s) = \frac{\uppi^{(d_1+d_2)/2}}
{\Gamma(\nu)\alpha_1} \biggl(\frac{1}{2}s\biggr)^{2\theta'} H^{2,1}_{2,2}
\left( \frac{(({1}/{2})s)^{2/\alpha_1}}{(({1}/{2})r)^{2}
}\left|
\matrix{(1,1), \biggl(\displaystyle\frac{d_1}{2},1\biggr) \vspace*{2pt}\cr\biggl(\displaystyle\frac{d_1}{2\alpha_1},\frac
{1}{\alpha_1}\biggr),\biggl(-\theta',\displaystyle\frac{1}{\alpha_1}\biggr)}\right.\right),
\label{foxh}
\end{equation}
where this result holds for $r=0$ by continuity.
For the parameter values
of the $H$-function in~(\ref{foxh}), from (\ref{astar}) and (\ref{Delta})
in Section 3.2,
$\Delta<0$ and $a^*>0,$ so that by~\cite{kilbas}, Theorem 1.11,
$H_{2,2}^{2,1}(z)$ has an asymptotic expansion
as $z\to0$ for $|\arg z| < \frac{1}{2}a^\ast\uppi$.
To avoid complications, let $\W$ be the set of nonnegative integers,
assume $\theta\notin\W$ and that for all $\ell\in\W$,
$(\theta'-\ell)\alpha_1\notin\W$ and $(\theta-\ell)\alpha_1\notin\N$.
Then, as $z\to0$,
%
\begin{equation}
H_{2,2}^{2,1}\left(z\left|
\matrix{(1,1), \biggl(\displaystyle\frac{d_1}{2},1\biggr) \vspace*{2pt}\cr\biggl(\displaystyle\frac{d_1}{2\alpha_1},\frac
{1}{\alpha_1}\biggr),\biggl(-\theta',\displaystyle\frac{1}{\alpha_1}\biggr)}\right.\right)
\sim\sum_{\ell=0}^\infty\bigl\{h_{1\ell}^*z^{\ell\alpha_1+d_1/2}+
h_{2\ell}^*z^{(\ell-\theta')\alpha_1}\bigr\},
\label{asymp}
\end{equation}
where, using the duplication formula for~$\Gamma$,
%
\begin{eqnarray}\label{h1l}
h_{1\ell}^* & = & \frac{(-1)^\ell\alpha_1\Gamma(-\theta-\ell)
\Gamma(({1}/{2})d_1+\ell\alpha_1)} {\ell!\Gamma(-\ell\alpha_1)}
\nonumber
\\[-8pt]
\\[-8pt]
\nonumber
& = & \frac{\alpha_1\sin(\uppi\ell\alpha_1)\Gamma(({1}/{2})d_1+\ell
\alpha_1)\Gamma(\ell\alpha_1+1)}{\sin(\uppi\theta)\Gamma(\theta+\ell+1)
\ell!}
\end{eqnarray}
and
\begin{eqnarray*}
h_{2\ell}^*  &= & \frac{(-1)^\ell\alpha_1\Gamma(\theta-\ell) \Gamma
((\ell-
\theta')\alpha_1)} {\ell!\Gamma(\alpha_1(\theta-\ell))}\\
& =&  \frac{\alpha_1\sin\{\uppi\alpha_1(\theta-\ell)\}\Gamma((\ell-\theta')
\alpha_1)\Gamma((\ell-\theta)\alpha_1+1)}{\sin(\uppi\theta)\Gamma(1-\theta
+\ell)
\ell!}.
\end{eqnarray*}
From (\ref{h1l}), $h_{10}^\ast=0$ and, if $\alpha_1$ is an
integer, $h_{1\ell}^*=0$ for all $\ell$. For $\theta\in\W$,
\cite{kilbas}, (1.8.2)
applies, yielding an asymptotic expansion with
logarithmic terms, but the result is rather messy and is omitted here.

When $(\frac{1}{2}s)^{2/\alpha_1}/(\frac{1}{2}r)^2$
is small, (\ref{Gres}) converges
slowly and is numerically unstable.
Specifically, for any fixed $m>\alpha_1\theta'$ and $r>0$, as
$s\downarrow0$,
the $m$th term in (\ref{Gres}) tends to $\pm\infty$, even though $
\lim_{s\downarrow0}\tilde G(r,s)\to\tilde G(r,0)$, which is finite.
Thus, there must be a near canceling of large terms of opposite signs in
(\ref{Gres}) for $s$ small, so that high precision
arithmetic would be needed to obtain accurate results for $s$ sufficiently
small.
Fortunately,
the asymptotic expansion (\ref{asymp}) can be used to approximate
$\tilde G(r,s)$ for $s$ small.
Substituting (\ref{asymp}) into (\ref{foxh}) and considering $r>0$ fixed,
%
\begin{eqnarray}\label{Gasymp}
\tilde G(r,s)  \sim \frac{\uppi^{(d_1+d_2)/2}}{\Gamma(\nu)\alpha_1}
\sum_{\ell=0}^\infty\biggl\{ h_{1\ell}^*
\biggl(\frac{1}{2}s\biggr)^{2\theta}
\biggl(\frac{1}{2}r\biggr)^{-d_1} + h_{2\ell}^*
\biggl(\frac{1}{2}r\biggr)^{2\alpha_1\theta'} \biggr\}
\biggl\{\frac{(({1}/{2})s)^2}{(({1}/{2})r
)^{2\alpha_1}}
\biggr\}^\ell
\end{eqnarray}
as $s\downarrow0$.
Since $h_{10}^\ast=0$,
when $\theta$ is not an integer,
%
\begin{equation}
\tilde G(r,s) = \sum_{\ell=0}^{\lfloor\theta\rfloor+1}
A_\ell s^{2\ell} + B\gamma_{\theta+1}(s)+ \mathrm{o}(s^{2\theta})
\label{pit}
\end{equation}
as $s\downarrow0$
for some constants (depending on $r$)
$A_0,\ldots,A_{\lfloor\theta\rfloor+1},B$, where $B=0$ if and only
if $\alpha_1$ is an integer.
The asymptotic expansion (\ref{pit}) also holds when $\theta$ is an integer
by~\cite{kilbas}, (1.8.2) and~(1.4.5).

It is apparent that (\ref{Gres}) does not give a valid power series expansion
for $\tilde G$ when $\alpha_1<1$,
since it is easy to show that the individual terms in the
sum do not tend to 0 as $m\to\infty$ for fixed and positive $r$ and $s$.
Nevertheless, the representation in terms of $H$-functions
given by (\ref{foxh}) may still be valid for $\alpha_1\le1$
when $2\alpha_1\theta> d_1$ (so that $\theta' > 0$).
Excluding values of $(\alpha_1,\nu)$ for which (\ref{Hcondition}) is
not satisfied, a natural conjecture is that (\ref{foxh}) holds
when $\alpha_1< 1$ if $2\alpha_1\theta> d_1$.
A plausible approach to proving this conjecture
would be to use analytic continuation, but this would require at the least
extending the definition of $H$-functions to a strip of complex values of
$\alpha_1$ containing the positive real axis.
Note that
when $\alpha_1<1$, (\ref{Delta}) implies $\Delta>0$ and, by
\cite{kilbas}, Theorem 1.3, (\ref{asymp}) becomes a convergent power
series for $H_{2,2}^{2,1}$.
On the other hand, now (\ref{Hseries}) is no longer a convergent power series,
but, with $a^\ast$ defined as in (\ref{astar}),
it is a valid asymptotic expansion as $z\to\infty$ for $|\arg z| <
\frac{1}{2} a^\ast\uppi$ by~\cite{kilbas}, Theorem 1.7.

For $\alpha_1=1$,
(\ref{foxh}) can be directly verified when (\ref{Hcondition}) holds.
Specifically, for $\alpha_1=1$, using~\cite{kilbas}, Property 2.2,
the right-hand side of (\ref{foxh}) reduces to
%
\begin{equation}
\frac{\uppi^{(d_1+d_2)/2}}{\Gamma(\nu)}\biggl(\frac{1}{2}s\biggr)^{2\theta'}
H_{1,1}^{1,1}\left(\displaystyle\frac{s^2}{r^2}\left|\matrix{(1,1)\vspace*{2pt}\cr(-\theta',1)}
\right.\right).
\label{alp1}
\end{equation}
The parameter $\Delta$ defined in (\ref{Delta}) equals 0
and, thus, one can show
that the $H$-function has a convergent
power series given by~\cite{kilbas}, Theorem~1.3, when $s<r$
and by~\cite{kilbas}, Theorem~1.4, when $s>r$.
Consider $s<r$ and $\theta'\notin\W$.
Then straightforward calculations yield that (\ref{alp1}) equals
\begin{eqnarray*}
& & \frac{\uppi^{(d_1+d_2)/2}}{\Gamma(\nu)2^{2\theta'}} r^{2\theta'}
\sum_{\ell=0}^\infty\frac{(-1)^\ell\Gamma(-\theta'+\ell)}{\ell!}
\biggl(\frac{s}{r}\biggr)^{2\ell}
\\
& &\quad = -\frac{\uppi^{(d_1+d_2+2)/2}}{\sin(\uppi\theta')\Gamma(\theta'+1)
\Gamma(\nu)2^{2\theta'}}r^{2\theta'}
\sum_{\ell=0}^\infty{\theta' \choose\ell}\biggl( \frac{s}{r}
\biggr)^{2\ell},
\end{eqnarray*}
which, using the binomial series, equals (\ref{isotropic}).
A similar argument shows (\ref{alp1}) equals (\ref{isotropic}) for $s>r$,
and it additionally holds for $s=r$ by continuity.

\subsection{Consequences}\label{sec2.4}

One goal of this paper was to find a class of generalized covariance functions
that has a member satisfying (\ref{coord1}) and (\ref{coord2})
for all $C_1,C_2,\zeta_1$ and $\zeta_2$ positive.
In fact, the functions of the
form $G(b_1\bx,b_2\by)$ with $\alpha_1\ge1$ and $b_1$ and $b_2$ positive
satisfy this requirement.
To prove this, first suppose $\zeta_1\ne\zeta_2$ and,
without loss of generality, take $\zeta_1>\zeta_2$.
Theorem~\ref{Gtheorem} gives explicit expressions for positive constants
$D_1$ and $D_2$ such that $G(\bx,\ze)=
D_1\gamma_{\alpha_1\theta'}(\frac{1}{2}|\bx|)$ and
$G(\ze,\by)=
D_2\gamma_{\theta'}(\frac{1}{2}|\by|)$.
Setting $\theta' = \zeta_2$ and $\alpha_1 = \zeta_1/\zeta_2$
achieves the desired degree of smoothness in all directions.
From (\ref{Gres}), $\tilde G(0,b_2s) = D_2(\frac{1}{2}b_2
)^{\zeta_2}
\gamma_{\zeta_2}(s)$,
so by appropriate choice of $b_2$, one obtains $\tilde G(0,b_2s) =
C_2 \gamma_{\zeta_2}(s)$.
For $\zeta_1=\alpha_1 \theta'$ not an integer, by (\ref{G1k0.0}),
$\tilde G(b_1r,0) = D_1(\frac{1}{2}b_1)^{\zeta_1} \gamma_{\zeta_1}(r)$,
so by appropriate choice of $b_1$, one gets $\tilde G(b_1r,0) =
C_1 \gamma_{\zeta_1}(r)$.
When $\zeta_1 = \alpha_1\theta'$ is an integer, by (\ref{G1k0.1}),
$\tilde G(b_1r,0) = D_1(\frac{1}{2}b_1)^{2\zeta_1}
\gamma_{\zeta_1}(r)$ plus some
constant times $r^{2\zeta_1}$, so that $b_1$ can be chosen to make
$\tilde G(b_1r,0) = C_1 \gamma_{\zeta_1}(r)$ plus some constant times
$r^{2\zeta_1}$.
When $\zeta_1=\zeta_2$, set $\alpha_1=\alpha_2=1$ and $\zeta_1 =
\nu-\frac{1}{2}d_1-\frac{1}{2}d_2$.
By (\ref{isotropic}),
there exists $D>0$ such that $G(b_1\bx,b_2\by) =
D\gamma_{\zeta_1}(\sqrt{|b_1\bx|^2+|b_2\by|^2})$.
As before, one can clearly choose $b_1$ and $b_2$ so that (\ref
{coord1}) and
(\ref{coord2}) are satisfied as long as one ignores an even polynomial
of degree $2\zeta_1$ when $\zeta_1$ is an integer.
Since a GC-$k$ is only identified up to
even polynomials of degree at most $2k$,
it is fair to say that the class of generalized covariance functions
corresponding to (\ref{model2}) with $\alpha_2=1$ and $\alpha_1\ge1$
includes members satisfying
(\ref{coord1}) and (\ref{coord2}) for all $\zeta_1\ge\zeta_2,C_1$ and $C_2$.

Now consider in what sense members of
$G$ achieve the goal, identified in Section~\ref{sec1},
of being smoother away from the origin than they are at the origin.
If a function is not infinitely differentiable in any direction at the
origin but is infinitely differentiable everywhere but the origin,
then one might say without
controversy that such a function is smoother away from the
origin than at the origin.
Thus, when $\alpha_1\in\N$, the issue is settled.
But when $\alpha_1$ is not an integer,
$B$ in (\ref{pit}) is not~0, and
$G$ is not infinitely differentiable away from the origin.

One way to describe the smoothness of a function is by its pointwise
H\"older exponent $s$.
Consider a function $f$ from $\R^d$ to $\R$, $s>0$ and $\bx_0\in\R^d$.
Then $f\in C^s(\bx_0)$ if and only if there exists $\varepsilon>0$ and
a polynomial $P$ of degree less than $\lfloor s\rfloor$ and a constant $C$
such that $|f(\bx)-P(\bx-\bx_0)|\le C|\bx-\bx_0|^s$ for all
$\bx$ satisfying $|\bx-\bx_0|< \varepsilon$.
The H\"older exponent of $f$ at $\bx_0$, which I will denote by $\HE(\bx_0,f)$,
equals $\sup\{s\dvt s\in C^s(\bx_0)\}$.
Because the degree of smoothness at the origin of $G$ varies in different
directions, it will not suffice to compare the H\"older exponent at the origin
to the H\"older exponent elsewhere.
To avoid the kind of anomaly described in Section~\ref{sec1}, consider
the smoothness of $G$ in each direction separately.
Specifically, for vectors $\bz_0,\bz_1\in\R^{d_1+d_2}$
and (generalized) covariance function $K$,
consider the function
of $t\in\R$ given by $\bar K(t;\bz_0,\bz_1) = K(\bz_0+t\bz_1)$.
Define $\ze$ to be a vector of zeroes whose length is apparent from
context.
Then I claim that, in the present setting,
a useful notion of $K$ being smoother away from the origin
than at the origin is
%
\begin{equation}
\HE(0,\bar K(\cdot;\bz_0,\bz_1)) >
\HE(0,\bar K(\cdot;\ze,\bz_1))\qquad \mbox{for all } \bz_0\ne\ze, \bz_1\ne\ze.
\label{HE}
\end{equation}

To see why this definition might be appropriate here,
consider the following generalization of the example in the introduction.
Suppose $K$ is a covariance function or a generalized covariance function
of order 0 and,
for nonzero $\bz_0$ and $\bz_1$ and $0<\alpha<2$, $K(t\bz_1)
= C_0+C_1|t|^\alpha+\mathrm{o}(|t|^\alpha)$ and $K(\bz_0+t\bz_1) =
D_0+D_1|t|^\alpha+D_2 t+\mathrm{o}(|t|^\alpha)$ as $t\to0$ for some constants
$C_0,C_1,D_0,D_1$ and $D_2$ (possibly depending on $\bz_0$ and $\bz_1$) with
$C_1$ and $D_1$ nonzero.
It follows that $\HE(0,\bar K(\cdot;\ze,\bz_1)) =
\HE(0,\bar K(\cdot;\bz_0,\bz_1)) = \alpha$.
Furthermore, $\lim_{t\to0}\corr\{Z(t\bz_1)-Z(\ze),
Z(\bz_0+t\bz_1)-Z(\bz_0)\} = D_1/C_1\ne0$.
Now suppose the lack of
smoothness of $K$ at $\bz_0$ is localized in the
sense that there exists $\varepsilon>0$ such
that $\HE(0,\bar K(\cdot;\bz_0+\delta\bz_1,\bz_1))>\alpha$ for all
$0<|\delta|<\varepsilon$, which, as far as I am aware, holds for any space--time
covariance function that has been proposed in the literature.
This condition implies
$\lim_{t\to0}\corr\{Z(t\bz_1)-Z(\ze),
Z(\bz_0+\delta\bz_1+t\bz_1)-Z(\bz_0+\delta\bz_1)\} = 0$ for all
$0<|\delta|<\varepsilon$.
Thus, when $K$ is not smoother in the $\bz_1$ direction at $\bz_0$ than
it is at $\ze$, there is a ``discontinuity'' in correlations of
increments.
If, instead, $D_1=0$, which will be the case under~(\ref{HE}),
then this limiting correlation is 0 for all $\delta$
in a neighborhood of 0 including $\delta=0$, and no discontinuity occurs.

Define $\bz_j = (\bx_j,\by_j)$ for $j=0,1$
with $\bx_j\in\R^{d_1}$ and $\by_j\in\R^{d_2}$.
For $G$ as given in Theorem~\ref{Gtheorem}, let $\bar G(t;
\bz_0,\bz_1) = G(\bz_0+t\bz_1)$.
For $\bz_1\ne\ze$, by (\ref{Gres}),
$\HE(0,\bar G(\cdot;\ze,(\bx_1,\by_1))=2\theta'$ if $\bx_1\ne\ze$
and $\HE(0,\bar G(\cdot;\ze,(\ze,\by_1))=2\alpha_1\theta'$.
Now consider $\bz_0\ne\ze$.
If $\by_0\ne\ze$, then it is possible to show (\ref{Gres}) can be
differentiated termwise and $\HE(0,\bar G(\cdot;(\bx_0,\by_0),(\bx_1,\by_1)))
=\infty$.
Next, $\HE(0,\bar
G(\cdot;(\bx_0,\ze),(\ze,\by_1)))
=\infty$ by (\ref{G1k0.0}) or (\ref{G1k0.1}).
Finally, if $\bx_1\ne\ze$, then $\HE(0,\bar
G(\cdot;(\bx_0, \ze),(\bx_1,\by_1))) \ge2\theta+2$
by (\ref{pit}).
Because $\theta>\theta'$, in all cases
$\HE(0,\bar G(\cdot;\bz_0,\bz_1)) >
\HE(0,\bar G(\cdot;\ze, \bz_1))+2$ for all nonzero $\bz_0$ and
$\bz_1$, so (\ref{HE}) is more than satisfied.

It is not clear that satisfying (\ref{HE}), or even the stronger
condition met by $G$ here, will exclude all possible
``discontinuities'' or other anomalies in the covariance structure,
but it does avoid at least the type considered here.
It might be preferable to find covariance functions that satisfy~(\ref{coord1}) and~(\ref{coord2}) and are infinitely differentiable
away from the origin, but for $\zeta_1, \zeta_2$ and $\zeta_1/\zeta_2$ all
irrational, I am unaware of any generalized covariance functions
that satisfy all of these conditions.

Next, consider the problem of the dimple in fully symmetric
(even in both its arguments)
stationary space--time covariance functions described in~\cite{kent}.
A formal definition of the dimple is given in~\cite{kent}, but
the essential point is that the space--time covariance function $K(\bx,t)$
has a dimple in, say, the time lag $t$ if, for some fixed spatial lag
$\bx$,
$K(\bx,t)$ has a local minimum in $t$ at $t=0$.
This dimple implies that at the spatial lag $\bx$,
correlation is stronger with
both the near future and the near past than with the present,
and~\cite{kent} argues that such a lack of monotonicity in the
covariance structure will often be undesirable.
A dimple in the spatial lag can be defined similarly.

For a GC-$k$ with $k>0$, it is not clear what one should mean by a dimple,
but for $k=0$, the GC-0 that equals 0 at the origin is just minus the
semivariogram for the process: $\frac{1}{2}\Var\{Z(\bx,t)-Z(\ze,0)\}
=-G(\bx,t)$.
If the variogram is expected to increase as one
moves ``farther away'' in space--time, then
GC-0s with dimples should be avoided.

Assume $k_0=0$ or, equivalently, $\alpha_1\theta'<1$, so that $\tilde G$
is a GC-0.
First, if $\alpha_1=1$, $\tilde G(r,s)$ is a decreasing function of
$\sqrt{r^2+s^2}$, and there is no dimple, so assume $\alpha_1>1$.
In addition, assume~$\theta$ is not an integer so that (\ref{Gasymp})
holds.
To show that $\tilde G(r,s)$ does not have a dimple, it suffices to show
that for every $s\ge0$, $\tilde G(0,s) > \tilde G(r,s)$ for all $r$
sufficiently small
and, for every $r\ge0$, $\tilde G(r,0) > \tilde G(r,s)$ for all $s$
sufficiently small.
That $\tilde G(0,0)=0$ is greater than $\tilde G(r,0)$ and $\tilde G(0,s)$
for all positive $r$ and $s$ is immediate from (\ref{Gres}) and
(\ref{G1k0.0}).
For fixed $s>0$, from (\ref{Gres}),
\[
\tilde G(0,s) - \tilde G(r,s)
= \frac{c_1(\alpha_1)}{4\alpha_1 \Gamma(\nu)}r^2\gamma_{\theta'-1/\alpha_1}
\biggl(\frac{1}{2}s\biggr) +\mathrm{O}(r^4)
\]
as $r\downarrow0$.
It follows from $c_1(\alpha_1)>0$ and $1/\alpha_1 -\theta' > 0$ that
$\tilde G(0,s) > \tilde G(r,s)$ for all $r$ sufficiently small.
From (\ref{Gasymp}), for fixed $r>0$,
\[
\tilde G(r,0) - \tilde G(r,s) =
\frac{\uppi^{(d_1+d_2)/2}\Gamma(\theta-1)\Gamma((1-\theta')\alpha_1)}
{\Gamma(\nu)\Gamma(\alpha_1(\theta-1))}\biggl(\frac{1}{2}r\biggr)^{2\alpha_1
(\theta'-1)}\biggl(\frac{1}{2}s\biggr)^2 +\mathrm{o}(s^2)
\]
as $s\downarrow0$.
For $\theta>1$, $\tilde G(r,0) - \tilde G(r,s)$ is clearly positive for all
$s$ sufficiently small, so now consider $0<\theta<1$. In
this case, $\Gamma(\theta-1) < 0$,
so $\tilde G(r,0) > \tilde G(r,s)$ follows for all
$s$ sufficiently small if $\Gamma(\alpha_1(\theta-1))<0$,
which holds if $\alpha_1(\theta-1) \in(-2m-1,-2m)$ for some nonnegative
integer $m$.
This condition does not hold for all $\alpha_1$ and $\theta'$ for which
$\alpha_1 > 1$ and $\alpha_1\theta' < 1$, but it does always hold when
$\theta' \ge\frac{2}{4+d_1}$.
To prove this, it suffices to show $\alpha_1(\theta-1) > -1$,
which holds if $\theta' > 1 - (1+\frac{1}{2}d_1)/\alpha_1$.
The curves $\theta' = 1 - (1+\frac{1}{2}d_1)/\alpha_1$
and $\theta' = 1/\alpha_1$ intersect at $(\alpha_1,\theta') =
(2 + \frac{1}{2}d_2,2/(4+d_1))$,
from which it follows that $\alpha_1(\theta-1) >
-1$ holds for all $\theta' < 1/\alpha_1$
whenever $\theta' \ge\frac{2}{4+d_1}$.
This lower bound is $\frac{2}{5}$ for $d_1=1$ and is smaller for
larger $d_1$.
The lower bound of $\frac{2}{5}$ may not be too restrictive
in practice: Brownian motion has generalized covariance function
proportional to $\gamma_{1/2}$,
and processes less smooth than Brownian motion
are somewhat uncommon in applications.

\section{Discussion}\label{sec3}

For $G$ as defined by (\ref{step0}),
the series expansions (\ref{Gres}) and (\ref{Gasymp}) should, in
principle, allow fast
and accurate calculation of $G$, but there do not appear to be any
publicly available programs for computing $H$-functions and writing
general purpose code to carry out these calculations would require a
major effort.
In particular, preliminary investigations suggest that
considerable care needs to be taken to piece together the convergent power
series (\ref{Gres}) and the asymptotic expansion~(\ref{Gasymp}) to
obtain accurate approximations for all values of the argument of the
function.
Further work would also be
needed to handle those values for $(\alpha_1,\nu)$ for which one of the
expansions has a singularity or near singularity, including values of
$\alpha_1$ near 1.

Restricting $\alpha_2=1$ in (\ref{model2}) was essential to the derivation
of series expansions for the resulting generalized covariance functions.
Model (\ref{model2}) was, in turn, a simplification of (\ref{model1}), which
includes two range parameters.
Perhaps $H$-functions can be used to express, in at least some cases,
the (generalized) covariance
functions corresponding to these more general models.
However, even the richer class of covariance functions given by (\ref{model1})
is inadequate for modeling many natural processes.
In particular, these covariance
functions all satisfy $G(\bx,\by)=G(\bx,-\by)$ and, hence,
are all what~\cite{gneiting} calls fully symmetric.
Any process with a predominant direction of flow will not be fully symmetric,
so this constraint is often inappropriate.
Covariance functions that are not fully symmetric can
be generated from covariance functions
that are~\cite{stein2005,jun},
and these approaches can, in principle, be applied to the models
considered here.
An easy extension of this model is to allow for geometric anisotropies in
either $\bx$ or $\by$ by considering $\tilde G(|A\bx|,|B\by|)$ for any
$d_1\times d_1$ matrix $A$ and any $d_2\times d_2$ matrix~$B$.\looseness=1

As noted in the \hyperref[sec1]{Introduction}, there has been quite a lot of research
in recent years developing new classes of space--time covariance functions.
In many of these works, the focus has been on obtaining simple closed form
expressions for space--time covariance functions.
Having closed form expressions is certainly valuable in applying
the models, but it is critical that any such model provides
a good description of the spatial-temporal variations of the
process to which it is to be applied.
Comparing various models in a broad range of applications is one
important way to learn about which models will be of most use in
practice, but it is also important to consider the theoretical
properties of these models, such as their smoothness properties, both
at the origin and away from the origin, and the presence of dimples
or other possible anomalies.
Finding covariance function models for space--time processes that
allow for a different degree of smoothness in space and in time,
possess certain other desirable properties such as (\ref{f-cond})
and are accurately computable using series expansions is a major
challenge.
The results obtained here perhaps provide a first step to show
that it may not be necessary to sacrifice desired theoretical
properties of space--time models in order to gain computational tractability,
although admittedly quite a bit of work on numerical methods would
be needed before the generalized covariance functions proposed here could
be used routinely (or even not so routinely) in practice.

\begin{appendix}
\section*{Appendices}\label{app}

\subsection{\texorpdfstring{Proof of Theorem \protect\ref{Gtheorem}}{Proof of Theorem 1}}\label{appa1}
To prove (\ref{Gres}) for the process $Z$ with spectral
density (\ref{model2}), consider the process
$Z_1(x_1,\by) = Z((x_1,0,\ldots,0),\by)$, which has GC-$k_0$
$G_1(x_1,\by)=G((x_1,0,\ldots,0),\by)=\tilde G(|x_1|,|\by|)$.
Assume for now that $k_0 = \lfloor\alpha_1\theta'
\rfloor$ does not equal $\alpha_1\theta'$
so that $\alpha_1\theta' - 1 < k_0 <
\alpha_1\theta' $.
Then the process $Z_1$
is $k_0$ times mean square differentiable in its first coordinate
direction and, for $m\le k_0$, denote its $m$th mean square derivative
process by $Z_1^m(x_1,\by)$.
The generalized covariance function for $Z_1^m(x_1,\by)$, denoted by
$G_1^m(x_1,\by)$, can be chosen to satisfy\vspace*{1pt}
\renewcommand{\theequation}{\arabic{equation}}
\setcounter{equation}{23}
\begin{equation}
G_1^{m}(x_1,\by)=
(-1)^{m} \frac{\partial^{2m}}{\partial x_1^{2m}}
G_1(x_1,\by).
\label{G1m}\vspace*{1pt}
\end{equation}
Now
$\frac{\partial^{2m+1}}{\partial x_1^{2m+1}}
G_1(0,\by)=0$ for $m < k_0$, so if one knew
$\frac{\partial^{2m}}{\partial x_1^{2m}}
G_1(0,\by)$ for $m\le k_0$, then $G_1(x_1,\by)$
could be recovered from $G_1^{k_0}(x_1,\by)$ by integration.
Then, since $G(\bx,\by)$ has a version that only depends on $\bx$ through
$|\bx|$, one can obtain $G$.

Suppose $d_1>1$.
The case $d_1=1$ requires a slightly different but easier argument.
The process $Z_1^{k_0}(x_1,\by)$ is an IRF-0,
so its GC-0 can be taken to equal negative the semivariogram of the process.
Denoting $|\by|$ by $s$,
\begin{eqnarray*}
 -G_1^{k_0}(x_1,\by)
&  = &\int_{\R^{d_1}}\int_{\R^{d_2}}
(1-\mathrm{e}^{\mathrm{i}\tau_1 x_1 + \mathrm{i}\bom'\by})\tau_1^{2k_0}
(|\btau|^{2\alpha_1}+|\bom|^2)^{-\nu}\bdom \bdtau\\
& = &\frac{\uppi^{d_2/2}}{2^{\theta-1}\Gamma(\nu)}
\int_{\R^{d_1}} \frac{\tau_1^{2k_0}}{|\btau|^{2\alpha_1\theta}}
\{ \mathcal{M}_\theta(0)-\cos(\tau_1 x_1)\mathcal{M}_\theta(|\btau
|^{\alpha_1}s)
\} \bdtau.\
\end{eqnarray*}
Switching to hyperspherical coordinates with $|\btau|=u$, $\tau_1=u\cos
\phi_1$
and integrating over the angles $\phi_2,\ldots,\phi_{d_1-1}$ yields
%
\begin{eqnarray}\label{G1k0}
 -G_1^{k_0}(x_1,\by)
&= &\frac{\uppi^{(d_1+d_2-1)/2}}{2^{\theta-2}\Gamma(\nu)\Gamma({(d_1-1)}
/{2}) }\nonumber\\
&&{}\times\int_0^\infty u^{2k_0-2\alpha_1\theta'-1}
\int_0^\uppi\cos^{2k_0} \phi_1 \sin^{d_1-2} \phi_1
\\
& &\hspace*{98pt} {}   \times
\{\mathcal{M}_\theta(0)-\cos(x_1 u\cos\phi_1)\mathcal{M}_\theta(u^{\alpha_1}s)
\} \,\mathrm{d}\phi_1  \,\mathrm{d}u.\nonumber
\end{eqnarray}
Making the change of variables $\sigma=\cos\phi_1$ and
using~\cite{gradshteyn}, 3.251.1 and 3.771.4,
and the series expansion for the generalized
hypergeometric function $_1  F_2$ yields
%
\begin{eqnarray}\label{step2}
& & \int_0^\uppi\cos^{2k_0} \phi_1 \sin^{d_1-2} \phi_1
\{\mathcal{M}_\theta(0)-\cos(x_1 u\cos\phi_1)\mathcal{M}_\theta(u^{\alpha_1}s)
\} \,\mathrm{d}\phi_1
\nonumber\\
& &\quad =  2\int_0^1 \sigma^{2k_0}(1-\sigma^2)^{(d_1-3)/2}\{\mathcal{M}_\theta(0)-\cos(x_1 u\sigma)\mathcal{M}_\theta(u^{\alpha_1}s)
\} \,\mathrm{d}\sigma\nonumber\\
& &\quad =   B\biggl(k_0+\frac{1}{2},\frac{d_1-1}{2}\biggr)
\nonumber
\\[-8pt]
\\[-8pt]
\nonumber
& &\qquad{}\times   \biggl\{
\mathcal{M}_\theta(0)-   _1  F_2\biggl(k_0+\frac{1}{2};\frac{1}{2},
k_0+\frac{d_1}{2};-\frac{x_1^2u^2}{4}\biggr)\mathcal{M}_\theta(u^{\alpha
_1} s)
\biggr\}
\\
& &\quad =   B\biggl(k_0+\frac{1}{2},\frac{d_1-1}{2}\biggr)
\Biggl[\{ \mathcal{M}_\theta(0)-\mathcal{M}_\theta(u^{\alpha_1} s)\}
\nonumber\\
& &\phantom{\quad =   B\biggl(k_0+\frac{1}{2},\frac{d_1-1}{2}\biggr)
\Biggl[}{}-\sum_{\ell=1}^\infty\frac{(k_0+{1}/{2})_\ell}
{(k_0+({1}/{2})d_1)_\ell(2\ell)!}( -x_1^2 u^2)^\ell
\mathcal{M}_\theta(u^{\alpha_1} s)\Biggr],\nonumber
\end{eqnarray}
where $B$ is the beta function.

The following properties of $\mathcal{M}_\theta$ are used in the proof.
For any $\theta>0$, as $t\downarrow0$,
%
\begin{equation}
\mathcal{M}_\theta(t)=\sum_{r=0}^{\lfloor\theta\rfloor} U_r
t^{2r}+V\gamma_\theta(t)+\mathrm{o}(t^{2\theta})
\label{matern.series}
\end{equation}
for appropriate values of the $U_r$s and $V$~\cite{stein-book}, Section 2.7.
In addition, for any $\theta>0$, there exist positive constants $C$ and $D$
(depending on $\theta$) such that
%
\begin{equation}
0\le\mathcal{M}_\theta(t)\le C\mathrm{e}^{-Dt}
\label{matern.bound}
\end{equation}
for all $t\ge0$.
Furthermore, for all real $\theta$ and all $t>0$,
%
\begin{equation}
\frac{\mathrm{d}}{\mathrm{d}t}\mathcal{M}_\theta(t) = -t\mathcal{M}_{\theta-1}(t).
\label{matern.deriv}
\end{equation}

Assume $s>0$ for now.
The function
$u^{2k_0-2\alpha_1\theta'-1}
\{ \mathcal{M}_\theta(0)-\mathcal{M}_\theta(u^{\alpha_1} s)\}$
is integrable in $u$ over $(0,\infty)$ since, from (\ref{matern.bound}) and
$k_0<\alpha_1\theta'$, it is integrable over $(0,1]$, and, from
(\ref{matern.series}) and $k_0>\alpha_1\theta'-1$, it is integrable
over $(1,\infty)$.
Furthermore,\vspace*{1.5pt}
\[
\sum_{\ell=1}^\infty\frac{(k_0+{1}/{2})_\ell}
{(k_0+({1}/{2})d_1)_\ell(2\ell)!}( x_1^2 u^2)^\ell
\le\cosh(x_1 u)\vspace*{1.5pt}
\]
and $\cosh(x_1 u) \mathcal{M}_\theta(u^{\alpha_1} s)$ is integrable over
$(0,\infty)$ for $\alpha_1>1$, so that\vspace*{1.5pt}
\[
\sum_{\ell=1}^\infty\frac{(k_0+{1}/{2})_\ell}
{(k_0+({1}/{2})d_1)_\ell(2\ell)!}( -x_1^2 u^2)^\ell
\mathcal{M}_\theta(u^{\alpha_1} s)\vspace*{1.5pt}
\]
can be integrated termwise over $(0,\infty)$.
Making the change of variables $t=u^{\alpha_1}$, integrating
by parts and using the definition of $\mathcal{M}_\theta$ gives\vspace*{1.5pt}
\begin{eqnarray*}
& & \int_0^\infty u^{2k_0-2\alpha_1\theta'-1} \{ \mathcal{M}_\theta(0)-\mathcal{M}_\theta(u^{\alpha_1} s)\} \,\mathrm{d}u
\\[1.5pt]
& & \quad= \frac{1}{\alpha_1}\int_0^\infty t^{2k_0/\alpha_1-\theta'-1}
\{ \mathcal{M}_\theta(0)-\mathcal{M}_\theta(ts)\} \,\mathrm{d}t
\\[1.5pt]
& &\quad = -\frac{s^2}{2\alpha_1\theta' - 2k_0}
\int_0^\infty t^{2k_0/\alpha_1 - 2\theta' + 1}\mathcal{M}_{\theta-1}(ts)\,\mathrm{d}t
\\[1.5pt]
& & \quad= -\frac{s^{\theta+1}}{2\alpha_1\theta' - 2k_0}
\int_0^\infty t^{(2k_0+d_1)/\alpha_1 - \theta}\mathcal{K}_{\theta-1}(ts)\,\mathrm{d}t.\vspace*{1.5pt}
\end{eqnarray*}
Since $(2k_0+d_1)/\alpha_1 - \theta> |1-\theta|$,
\cite{gradshteyn}, 6.561.16, applies, and one obtains\vspace*{1.5pt}
%
\begin{eqnarray}\label{term1}
& & \int_0^\infty u^{2k_0-2\alpha_1\theta'-1} \{ \mathcal{M}_\theta(0)-\mathcal{M}_\theta(u^{\alpha_1} s)\} \,\mathrm{d}u
\nonumber
\\[-4pt]
\\[-4pt]
\nonumber
& &\quad  = -\frac{2^{\theta-2}}{
\alpha_1}\Gamma\biggl( \frac{k_0}{\alpha_1}
-\theta'\biggr) \Gamma\biggl( \frac{2k_0+d_1}{2\alpha_1}\biggr)
\biggl(\frac{1}{2}s\biggr)^{2\theta'-2k_0/\alpha_1}.\vspace*{1.5pt}
\end{eqnarray}
Next, for $\ell\ge1$, again using the change of variables $t=u^{\alpha_1}$
and~\cite{gradshteyn}, 6.561.16, gives\vspace*{1.5pt}
%
\begin{eqnarray}\label{term2}
& & \int_0^\infty u^{2k_0+2\ell-2\alpha_1\theta'-1}
\mathcal{M}_\theta(u^{\alpha_1} s) \,\mathrm{d}u
\nonumber
\\[-4pt]
\\[-4pt]
\nonumber
& & \quad= \frac{2^{\theta-2}}{
\alpha_1}\Gamma\biggl( \frac{k_0+\ell}{\alpha_1}
-\theta'\biggr)  \Gamma\biggl( \frac{2k_0+2\ell+d_1}{2\alpha_1}\biggr)
\biggl(\frac{1}{2}s\biggr)^{2\theta'-(2k_0+2\ell)/\alpha_1}.\
\end{eqnarray}
For $\by\ne\ze$, using (\ref{G1k0}), (\ref{step2}), (\ref{term1})
and (\ref{term2}),
%
\begin{eqnarray}\label{G1k0.res}
& & G_1^{k_0}(x_1,\by)
\nonumber\\
& &\quad =
\frac{\uppi^{d_1+d_2-1}\Gamma((k_0+{1}/{2})}
{\Gamma(\nu)\Gamma(k_0+({1}/{2})d_1)\alpha_1}
\nonumber\\
& & \qquad{} \times\sum_{\ell=0}^\infty\frac{(k_0+{1}/{2}
)_\ell}
{(k_0+({1}/{2})d_1)_\ell(2\ell)!}
\Gamma\biggl(\frac{2k_0+2\ell+d_1}{2\alpha_1}\biggr) (-x_1^2)^\ell\gamma_{\theta'-(k_0+\ell)/\alpha_1}
\biggl(\frac{1}{2}s\biggr)
\\
& & \quad= \frac{1}{\alpha_1\Gamma(\nu)}\sum_{\ell=0}^\infty
\frac{(k_0+\ell)!({1}/{2})_{k_0+\ell}
c_{k_0+\ell}(\alpha_1)}{(2\ell)!} (-x_1^2)^\ell
\gamma_{\theta'-(k_0+\ell)/\alpha_1}
\biggl(\frac{1}{2}s\biggr)
\nonumber\\
& & \quad= \frac{1}{\alpha_1\Gamma(\nu)}\sum_{m=k_0}^\infty
\frac{m!({1}/{2})_m
c_m(\alpha_1)}{\{2(m-k_0)\}!} (-x_1^2)^\ell
\gamma_{\theta'-m/\alpha_1}
\biggl(\frac{1}{2}s\biggr).
\nonumber
\end{eqnarray}

Next consider $G_1^m(\ze,\by)$ for $m < k_0$.
The process $\frac{\partial^m}{\partial x_1^m}Z(\bx,\by)$ has spectral density
$\tau_1^{2m}(|\btau|^{2\alpha_1}+|\bom|^2)^{-\nu}$,
and hence the process $\frac{\partial^m}{\partial x_1^m}Z(\ze,\by)
=Z_1^m(0,\by)$ considered just
as a function of $\by\in\R^{d_2}$ has spectral density
%
\begin{eqnarray}\label{hyper}
& & \int_{\R^{d_1}} \frac{\tau_1^{2m}}{(|\btau|^{2\alpha_1}+|\bom
|^2)^\nu}
\bdtau\nonumber\\
& &\quad = \frac{2\uppi^{(d_1-1)/2}}{\Gamma({(d_1-1)}/{2})}
\int_0^\infty\int_0^\uppi\frac{u^{2m+d_1-1}\cos^{2m}\phi_1 \sin^{d_1-2}
\phi_1}{(u^{2\alpha_1}+|\bom|^2)^\nu} \,\mathrm{d}\phi_1  \,\mathrm{d}u
\\
& & \quad= \frac{\uppi^{(d_1-1)/2}\Gamma(m+{1}/{2}) B({(2m+d_1)}/
{2\alpha_1},\nu- {(2m+d_1)}/ {2\alpha_1})}{\alpha_1\Gamma(m+
({1}/{2})d_1) |\bom|^{2\nu-(2m+d_1)/\alpha_1}}
\nonumber
\end{eqnarray}
by switching $\btau$ to
hyperspherical coordinates and using~\cite{gradshteyn}, 3.241.4.
The process $Z_1^m(0,\by)$ with spectral density (\ref{hyper})
is an IRF-$(k_0-m)$ (it may be an IRF of lower
order as well) and its corresponding GC-$(k_0-m)$ can be taken as
\cite{gelfand}, Chapter II, Section 3.3, equations (2) and (11),
%
\begin{eqnarray}\label{boundary}
G_1^m(0,\by) & = & \frac{c_m(\alpha_1) m!
\Gamma(m+{1}/{2})}
{\uppi^{1/2}\alpha_1\Gamma(\nu)}
\gamma_{\theta'-m/\alpha_1}\biggl(\frac{1}{2}|\by|\biggr)
\nonumber
\\[-8pt]
\\[-8pt]
\nonumber
& = & \frac{c_m(\alpha_1) m!
({1}/{2})_m}
{\alpha_1\Gamma(\nu)}
\gamma_{\theta'-m/\alpha_1}\biggl(\frac{1}{2}|\by|\biggr).
\end{eqnarray}

To recover $G_1(x_1,\by)$ and hence $G(\bx,\by)$, repeatedly integrate
(\ref{G1k0.res}) and use (\ref{boundary}) to set the boundary conditions.
Specifically, $G_1(x_1,\by)$ must be of the form
\[
G_1(x_1,\by)  =  \sum_{\ell=0}^{k_0-1} x_1^{2\ell} F_\ell(\by) +(-1)^{k_0}\int_0^{x_1}\int_0^{z_1}\cdots
\int_0^{z_{2k_0-1}} G_1^{k_0}(z_{2k_0},\by)\,\mathrm{d}z_{2k_0}\cdots \,\mathrm{d}z_1
\]
for some suitable functions $F_0,\ldots,F_{k_0-1}$.
Substituting the series for $G_1^{k_0}$ in (\ref{G1k0.res}) into the preceding
expression and integrating termwise,
which is easily justified for $s>0$ by dominated convergence, yields
%
\begin{eqnarray}\label{G1.1}
 G_1(x_1,\by) &= &\sum_{\ell=0}^{k_0-1} x_1^{2\ell} F_\ell(\by)
\nonumber
\\[-8pt]
\\[-8pt]
\nonumber
& &{}
+ \frac{1}{\alpha_1\Gamma(\nu) }
\sum_{m=k_0}^\infty\frac{c_m(\alpha_1)m!({1}/{2} )_m
2^{2m}} {(2m)!}
\biggl\{-\biggl(\frac{1}{2}x_1\biggr)^2\biggr\}^m
\gamma_{\theta'-m/\alpha_1}\biggl(\frac{1}{2}|\by|\biggr).
\end{eqnarray}
Elementary calculations demonstrate
%
\begin{equation}
\frac{m!({1}/{2})_m 2^{2m}}{(2m)!}=1
\label{comb.id}
\end{equation}
for all $m\in\W$.
For $0\le m < k_0$, differentiating (\ref{G1.1}) $2m$ times,
setting $\by=\ze$ and using (\ref{boundary}) and~(\ref{comb.id}) gives
%
\begin{eqnarray}\label{Fm}
F_m(\by) & = & \frac{(-1)^m}{(2m)!}G_1^m(0,\by)
\nonumber
\\[-8pt]
\\[-8pt]
\nonumber
& = & \frac{(-{1}/{4})^m c_m(\alpha_1)}{\alpha_1\Gamma(\nu)}
\gamma_{\theta'-m/\alpha_1}\biggl(\frac{1}{2}|\by|\biggr).
\end{eqnarray}
Substituting (\ref{comb.id}) and (\ref{Fm}) into (\ref{G1.1})
yields
\[
G_1(x_1,\by) = \frac{1}{\alpha_1\Gamma(\nu)}
\sum_{m=0}^\infty \biggl\{-\biggl(\frac{1}{2}x_1\biggr)\biggr\}^m
c_m(\alpha_1)\gamma_{\theta'-m/\alpha_1}\biggl(\frac{1}{2}|\by|\biggr),
\]
and (\ref{Gres}) follows from $G_1(x_1,\by)=\tilde G(|x_1|,|\by|)$.

To obtain an explicit expression for $\tilde G(r,0)$,
go back to (\ref{G1k0}) and change the order of
integration.
Integrating by parts and using~\cite{gradshteyn}, 3.761.4,
\begin{eqnarray*}
& &
\int_0^\infty u^{2k_0-2\alpha_1\theta'-1}\{1-\cos(x_1 u\cos\phi_1)\} \,\mathrm{d}u
\\
& &\quad = \frac{x_1\cos\phi_1}{2\alpha_1\theta'-2k_0}\int_0^\infty
u^{2k_0-2\alpha_1\theta'}\sin(x_1 u \cos\phi_1) \,\mathrm{d}u \\
& & \quad= -\Gamma(2k_0-2\alpha_1\theta')
\sin\{\uppi(k_0-\alpha_1\theta')\} |x_1\cos\phi_1|^{
2\alpha_1\theta'-2k_0}.
\end{eqnarray*}
When $k_0\ne\alpha_1\theta'$,
substituting this result and $\mathcal{M}_\theta(0) = 2^{\theta-1}\Gamma
(\theta)$
into (\ref{G1k0}) and integrating over~$\phi_1$ yields
\begin{eqnarray*}
G_1^{k_0}(x_1,\ze) & = & \frac{2\uppi^{(d_1+d_2-1)/2}{\Gamma(\theta)
\Gamma(\alpha_1\theta'+{1}/{2})\Gamma(2k_0-2\alpha_1\theta')
}}
{\Gamma(\nu) \Gamma(\alpha_1\theta)}\\
&&{}\times \cos(\uppi\alpha_1\theta') (-1)^{k_0}|x_1|^{2\alpha_1\theta'-2k_0}.
\end{eqnarray*}
Now, $\Gamma(2k_0-2\alpha_1\theta')=(2\alpha_1\theta'-2k_0+1)_{2k_0}
\Gamma(-2\alpha_1\theta')$, so
applying the duplication formula for~$\Gamma$ to $\Gamma(-2\alpha
_1\theta')$
and then the reflection formula to $\Gamma(-\alpha_1\theta'+\frac
{1}{2})$
yields
\begin{eqnarray*}
&&G_1^{k_0}(x_1,\ze) \\
&&\quad    =\frac{(-1)^{k_0}\uppi^{(d_1+d_2)/2}\Gamma(\theta)
\Gamma(-\alpha_1\theta')(2\alpha_1\theta'-2k_0+1)_{2k_0}}
{\Gamma(\nu)\Gamma(\alpha_1\theta)2^{\alpha_1\theta'}}
|x_1|^{2\alpha_1\theta'-2k_0}.
\end{eqnarray*}
Integrating this expression $2k_0$ times and using the boundary condition
$G_1^m(0,\ze)=0$ for $m<k_0$ to make $G_1^m(0,\by)$ continuous at
$\by=\ze$ gives (\ref{G1k0.0}).

To show that (\ref{Gres}) holds when $k_0 = \alpha_1\theta'$,
write $\tilde G_{\nu}$ to make the dependence of $\tilde G$
on $\nu$ explicit (but still suppressing the dependence on
$\alpha_1$, $d_1$ and $d_2$).
Define $\nu_0= (k_0
+\frac{1}{2}d_1)/\alpha_1 + \frac{1}{2}d_2$ and view $\alpha_1$ as fixed.
For any given ALC-$k_0$,
the last line of (\ref{step0}) is continuous as $\nu\downarrow\nu_0$,
so the first line is as well.
Thus,
%
\begin{equation}
\sum_{\ell,j=1}^n \lambda_\ell\lambda_j
\tilde G_{\nu_0}(r_{\ell j},s_{\ell j})
= \lim_{\nu\downarrow\nu_0} \sum_{\ell,j=1}^n
\lambda_\ell\lambda_j \tilde G_{\nu}
(r_{\ell j},s_{\ell j}).
\label{G.lim}
\end{equation}
Let $M_0$ be the set of nonnegative integers $m$ for which $h_m
=(k_0-m)/\alpha_1 \in\W$.
This set is finite and includes $k_0$ as its largest element.
Writing $\nu= \nu_0 + \varepsilon$, define
\[
P_\varepsilon(r,s) = \sum_{m\in M_0} \frac{c_m(\alpha_1)\Gamma
(-h_m-\varepsilon)}
{\alpha_1\Gamma(\nu_0+\varepsilon)} \biggl\{-\biggl(\frac{1}{2}r\biggr)^2\biggr\}^m
\biggl(\frac{1}{2}s\biggr)^{2h_m}.
\]
Because $\sum_{\ell=1}^n\lambda_\ell Z(\bx_\ell)$ is an ALC-$k_0$,
for $a,b\in\W$ and $a+b\le k_0$,
subtracting any linear combination of terms\vspace*{1pt} like
$r_{\ell j}^{2a}s_{\ell j}^{2b}$ from $\tilde G_\nu(r_{\ell j},s_{\ell j})$
on the right-hand side of (\ref{G.lim}) does not change the result.
Now $\alpha_1>1$ implies $m+h_m\le k_0$, so $P_\varepsilon(r_{\ell
j},s_{\ell j})$
is of this required form and
\[
\sum_{\ell,j=1}^n \lambda_\ell\lambda_j
\tilde G_{\nu}(r_{\ell j},s_{\ell j}) =
\sum_{\ell,j=1}^n
\lambda_\ell\lambda_j \{\tilde G_{\nu}
(r_{\ell j},s_{\ell j}) - P_\varepsilon(r_{\ell j},s_{\ell j})
\}.
\]
Thus, to prove (\ref{Gres}), it suffices to show that for all $r,s$
nonnegative,
%
\begin{equation}
\tilde G_{\nu_0}(r,s)
= \lim_{\nu\downarrow\nu_0}
\{\tilde G_{\nu} (r,s)-P_\varepsilon(r,s)\}.
\label{G.lim.0}
\end{equation}
For $s>0$, for all $\varepsilon$ sufficiently small,
\begin{eqnarray*}
\tilde G_{\nu} (r,s)-P_\varepsilon(r,s)
& = & \sum_{m=0}^\infty\frac{c_m(\alpha_1)}{\alpha_1\Gamma(\nu
_0+\varepsilon)}
\biggl\{-\biggl(\frac{1}{2}r\biggr)^2\biggr\}^{m} \\
& &\hspace*{16pt}{}  \times\Gamma(-h_m-\varepsilon)
\biggl(\frac{1}{2}s\biggr)^{2h_m}
\biggl[ \biggl(\frac{1}{2}s\biggr)^{2\varepsilon} - 1\{m\in M_0\}
\biggr].
\end{eqnarray*}
Dominated convergence justifies taking the limit
$\nu\downarrow\nu_0$ (equivalently, as $\varepsilon\downarrow0$)
inside this infinite sum.
For $m\notin M_0$, the limit is trivial, so consider $m\in M_0$.
By the reflection formula for~$\Gamma$, for $n\in\W$,
%
\begin{equation}
\Gamma(-n-\varepsilon) =
\frac{\uppi(-1)^{n+1}}{\sin(\uppi\varepsilon)
\Gamma(n+1+\varepsilon)}.
\label{reflect}
\end{equation}
Using this result and straightforward calculus yields
\[
\Gamma(-h_m-\varepsilon)\biggl(\frac{1}{2}s\biggr)^{2h_m}
\biggl\{ \biggl(\frac{1}{2}s\biggr)^{2\varepsilon} - 1\biggr\}
\to\gamma_{h_m}\biggl(\frac{1}{2}s\biggr)
\]
as $\nu\downarrow\nu_0$, establishing (\ref{G.lim.0}) when $s>0$.

It remains to establish (\ref{G.lim.0}) when $s=0$ with $\tilde
G_{\nu_0}(r,0)$ given by (\ref{G1k0.1}).
For all $\varepsilon$ sufficiently small,
applying (\ref{reflect}) to $\Gamma(-k_0-\alpha_1\varepsilon)$
and using $\Gamma(1-\varepsilon) = -\varepsilon\Gamma(-\varepsilon)$,
for all $\varepsilon$ sufficiently small,
%
\begin{eqnarray}\label{A1}
& & \tilde G_\nu(r,0) - P_\varepsilon(r,0) \nonumber\\
& &\quad   = \frac{(-1)^{k_0+1}\uppi^{(d_1+d_2+2)/2}
\Gamma(\theta_\nu)}
{\sin(\uppi\alpha_1\varepsilon)\Gamma(\nu)\Gamma(\alpha_1\theta_\nu)
\Gamma(k_0+1+\alpha_1\varepsilon)}
\biggl(\frac{1}{2}r\biggr)^{2k_0+2\alpha_1\varepsilon} \nonumber\\
& & \qquad{}-
\frac{c_{k_0}(\alpha_1)\Gamma(-\varepsilon)}{\alpha_1\Gamma(\nu)}
\biggl\{-\biggl(\frac{1}{2}r\biggr)^2\biggr\}^{k_0}
\nonumber
\\[-8pt]
\\[-8pt]
\nonumber
& &\quad   = -\frac{\uppi^{(d_1+d_2)/2}}{\Gamma(\nu)}
\biggl\{-\biggl(\frac{1}{2}r\biggr)^2\biggr\}^{k_0}\\
&&\hspace*{6pt}\qquad{}\times\biggl\{
\frac{\uppi\Gamma(\theta_\nu)}
{\sin(\uppi\alpha_1\varepsilon)
\Gamma(\alpha_1\theta_\nu)\Gamma(k_0+1+\alpha_1\varepsilon)}
\biggl[\biggl\{\biggl(\frac{1}{2}r\biggr)^{2\alpha_1\varepsilon}-1\biggr\}+1\biggr]
\nonumber\\
& &\hspace*{24pt}\qquad {} - \frac{
\Gamma({(d_1+2k_0)}/{2\alpha_1})
\Gamma(1-\varepsilon)} {\alpha_1\varepsilon\Gamma(k_0+({1}/{2})d_1)
k_0!}\biggr\}.\nonumber
\end{eqnarray}
As $\varepsilon\downarrow0$,
%
\begin{equation}
\frac{\uppi}{\sin(\uppi\alpha_1\varepsilon)}
\biggl\{\biggl(\frac{1}{2}r\biggr)^{2\alpha_1\varepsilon}-1\biggr\}\to
2\log\biggl(\frac{1}{2}r\biggr).
\label{A2}
\end{equation}
By the definition of the digamma function $\psi$, $\Gamma(x+\varepsilon) =
\Gamma(x)\{1+\varepsilon\psi(x)+\mathrm{O}(\varepsilon^2)\}$
as $\varepsilon\to0$ as long as $-x\notin\W$.
Then
%
\begin{eqnarray}\label{A3}
& & \frac{\uppi\Gamma(\theta_\nu)}
{\sin(\uppi\alpha_1\varepsilon)
\Gamma(\alpha_1\theta_\nu)\Gamma(k_0+1+\alpha_1\varepsilon)}
- \frac{
\Gamma({(d_1+2k_0)}/{(2\alpha_1)})
\Gamma(1-\varepsilon)}
{\alpha_1\varepsilon\Gamma(k_0+({1}/{2})d_1)k_0!}
\nonumber
\\
& & \quad  =
\frac{\Gamma({(d_1+2k_0)}/{(2\alpha_1)})}
{k_0!\Gamma(k_0+({1}/{2})d_1)}\nonumber\\
&&\qquad{}\times
\biggl[\frac{1+\varepsilon\psi({(d_1+2k_0)}/{(2\alpha_1)})}
{\alpha_1\varepsilon\{ 1+\alpha_1\varepsilon\psi(k_0+
({1}/{2})d_1)
\}
\{ 1+\alpha_1\varepsilon\psi(k_0+1)\}} -
\frac{1-\varepsilon\psi(1)}{\alpha_1\varepsilon}\biggr] +\mathrm{O}(\varepsilon)
\qquad\\
& &  \quad =
\frac{\Gamma({(d_1+2k_0)}/{(2\alpha_1)})}
{k_0!\Gamma(k_0+({1}/{2})d_1)}\nonumber\\
&&\qquad{}\times \biggl\{\frac{1}{\alpha_1}\psi\biggl(\frac{2k_0+d_1}{2\alpha_1}\biggr)
+\frac{1}{\alpha_1}\psi(1) - \psi\biggl(k_0+\frac{1}{2}d_1\biggr) -
\psi(k_0+1) \biggr\} + \mathrm{O}(\varepsilon).\nonumber
\end{eqnarray}
Thus, when $k_0=\alpha_1\theta'$ and $\tilde G_{\nu_0}(r,0)$ is
defined as in (\ref{G1k0.1}),
(\ref{A1})--(\ref{A3}) imply (\ref{G.lim.0}) holds for $s=0$.

\subsection{Properties of $H$-functions}\label{appa2}
This Appendix provides some background material on $H$-functions and is taken
from~\cite{kilbas}, Section 1.1.
Suppose $m,n,p$ and $q$ are integers satisfying $0\le m\le q, 0\le n\le p$,
$a_1,\ldots,a_p,b_1, \ldots b_q$ are complex numbers and
$\alpha_1,\ldots,\alpha_p,\beta_1,\ldots,\beta_q$ are positive reals.
Then the $H$-function $H_{p,q}^{m,n}$ is defined by, for complex $z$,
%
\begin{eqnarray}\label{H-transform}
& & H_{p,q}^{m,n} \left( z\left |
\matrix{(a_1,\alpha_1),\ldots,(a_p,\alpha_p) \vspace*{2pt}\cr(b_1,\beta_1),\ldots
,(b_q,\beta_q)}
\right.\right)
\nonumber
\\[-8pt]
\\[-8pt]
\nonumber
& &\quad =  \frac{1}{2\uppi \mathrm{i}}\int_\mathcal{L}
\frac{\prod_{j=1}^m\Gamma(b_j+\beta_j s)\prod_{i=1}^n\Gamma(1-a_j+\alpha
_j s)}
{\prod_{i=n+1}^p\Gamma(a_i+\alpha_j s)\prod_{j=m+1}^q\Gamma(1-b_j+\beta
_j s)}
z^{-s} \,\mathrm{d}s,
\end{eqnarray}
where~\cite{kilbas}, p.~2, gives the form of the contour $\mathcal{L}$,
and an empty product is defined to be 1.
For this integral to be well defined, none of the poles of the gamma
functions in the two products in the numerator of (\ref{H-transform}) may
coincide, or
%
\begin{equation}
\alpha_i(b_j+\ell)\ne\beta_j(a_i-k-1)
\label{poles}
\end{equation}
for $1\le i\le n,1\le j\le m$ and all $k,\ell\in\W$.
The validity of series expansions
of $H_{p,q}^{m,n}$ generally depends on the signs of
the following two quantities:
%
\begin{equation}
a^*=\sum_{i=1}^n \alpha_i-\sum_{i=n+1}^p \alpha_i + \sum_{j=1}^m
\beta_j-\sum_{j=m+1}^q \beta_j
\label{astar}
\end{equation}
and
%
\begin{equation}
\Delta=\sum_{j=1}^q \beta_j-\sum_{i=1}^p \alpha_i,
\label{Delta}
\end{equation}
where an empty sum is defined to be 0.
\end{appendix}

\section*{Acknowledgement}
This research was supported by US Department of Energy Grant DE-SC0002557.


%

\printhistory

\end{document}